\documentclass{amsart} 
\usepackage{epsf} 
\usepackage{amssymb,latexsym, amsmath, amscd, array, graphicx} 
\def\hepsffile{\leavevmode\epsffile} 

\swapnumbers 
\numberwithin{equation}{section} 

\theoremstyle{plain} 
\newtheorem{thm}{Theorem}[section]

\theoremstyle{definition} 
\newtheorem{defin}[thm]{Definition}

\newtheorem{rem}[thm]{Remark}

\def\sign{\protect\operatorname{sign}}

\def\rvu{\protect\operatorname{rvu}}
\def\sign{\protect\operatorname{sign}}
\def\vb{\protect\operatorname{vb}}
\def\vu{\protect\operatorname{vu}}

\def\R{{\mathbb R}} 
\def\N{{\mathbb N}}

\def\1{\hbox{\rm\rlap {1}\hskip.03in{\rom I}}} 
\def\Bbbone{{\rm1\mathchoice{\kern-0.25em}{\kern-0.25em} 
{\kern-0.2em}{\kern-0.2em}I}} 

\begin{document}
\title[Virtual Bridge Number One Knots] 
{Virtual Bridge Number One Knots} 
\author[E.~Byberi and V.~Chernov (Tchernov)]{Evarist Byberi and Vladimir V. 
Chernov 
(Tchernov)}
\address{E.~Byberi, HB 0552 Dartmouth College, Hanover NH 03755, 
USA} 
\email{Evarist.Byberi@dartmouth.edu} 

\address{V. Chernov, Department of Mathematics, 
6188 Bradley Hall, Dartmouth College, Hanover NH 03755, 
USA} 
\email{Vladimir.Chernov@dartmouth.edu} 

\begin{abstract}
We define the virtual bridge number $\vb(K)$ and the virtual unknotting number $\vu(K)$ invariants 
for virtual knots. For ordinary knots $K$ 
they are closely related to the bridge number $b(K)$ and the unknotting number $u(K)$ and  we have $\vu(K)
\leq u(K), \vb(K)\leq b(K).$

There are no ordinary knots $K$ with $b(K)=1.$ We show there are infinitely many homotopy classes of virtual 
knots each of which contains infinitely many 
isotopy classes of $K$ with $\vb(K)=1.$ In fact for each $i\in \N$ there exists $K$ virtually homotopic (but 
not virtually isotopic) 
to the unknot with $\vb(K)=1$ and $\vu(K)=i.$ 

\end{abstract}


\maketitle

\section{Introduction into Virtual Knots}
The Virtual Knot Theory was introduced by L.~Kauffman~\cite{Kauffman}. Let us recall some of its basic 
notions.

A knot is a smooth embedding $S^1\to \R^3$. It can be described by its knot diagram which 
is a generic immersion of $S^1$ into the $\R^2$-plane equipped with information
about over-passes and under-passes at double points. A knot diagram $D$ gives rise to 
a Gauss diagram $G_D$ that is a circle parameterizing the knot with each pair of preimages of double points 
of $D$ connected by an oriented signed chord. The chords are oriented from the preimage point on the 
over-passing branch to the preimage point on the under-passing branch.
The sign of a chord is the sign of the
corresponding double point. The resulting Gauss diagram $G_D$ 
is called {\em the Gauss diagram of the knot diagram $D$.\/}

\begin{figure}[htbp]
\begin{center}
\epsfxsize 10cm \hepsffile{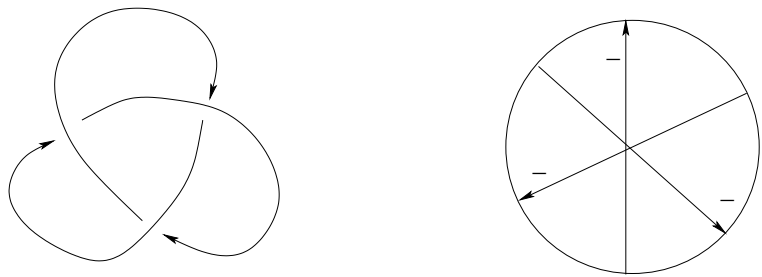}
\end{center}
\caption{A knot diagram and the corresponding Gauss diagram}\label{Example.fig}
\end{figure}

\begin{figure}[htbp]
\begin{center}
\epsfxsize 12cm \hepsffile{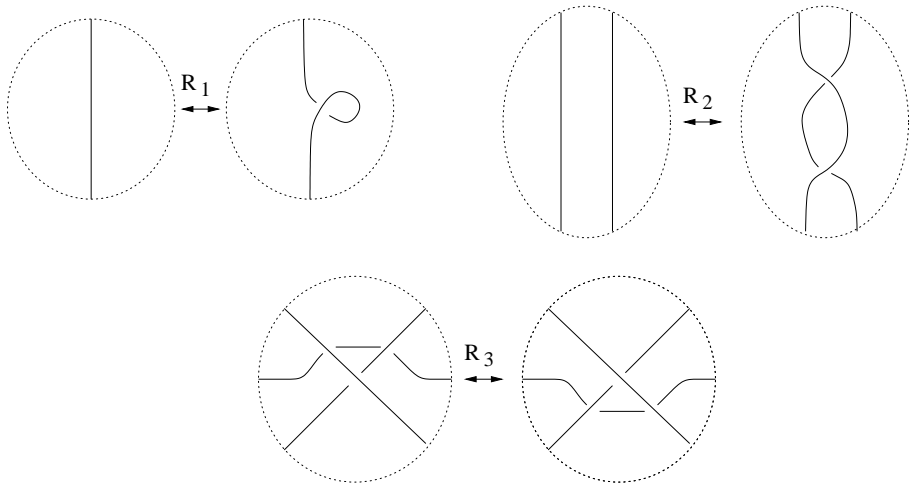}
\end{center}
\caption{Redemeister moves}\label{Redemeister.fig}
\end{figure}

\begin{figure}[htbp]
\begin{center}
\epsfxsize 8cm \hepsffile{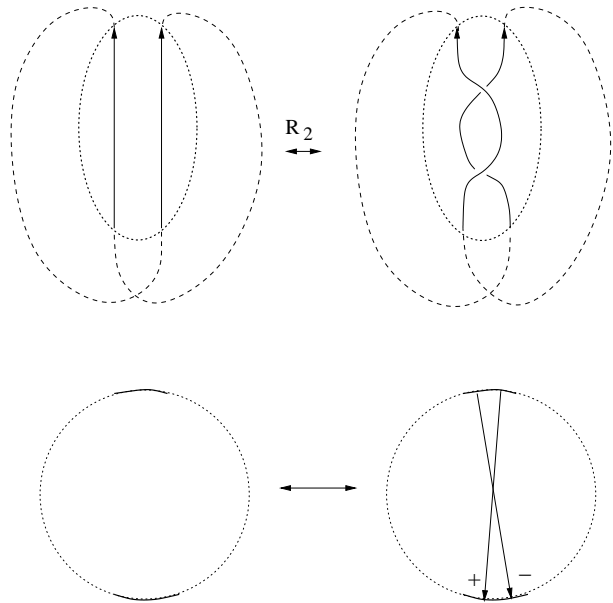}
\end{center}
\caption{Redemeister move and its encoding in terms of Gauss diagrams}\label{GaussRedemeister.fig}
\end{figure}

Gauss diagrams that are obtainable as Gauss diagrams of some knot diagrams
are said to be {\em realizable.\/} 
A knot diagram corresponding to a realizable Gauss diagram can be recovered
only up to a certain ambiguity. However the isotopy type of the corresponding
knot is recoverable in a unique way.

It is well-known that two knots described by their knot diagrams are isotopic
if and only if one can change one diagram to another by a sequence of ambient isotopies of the diagram and 
of {\em Reidemeister moves\/} shown in Figure~\ref{Redemeister.fig}. 

The Redemeister moves can be easily encoded~\cite{GPV} in the language of Gauss diagrams, see 
Figure~\ref{GaussRedemeister.fig}. The equivalence classes of the, not necessarily realizable, 
Gauss diagrams modulo the resulting moves are called virtual knots.

\section{Main results}
\begin{defin}[virtual bridge number]
Recall that a bridge in a knot diagram is an arc between two consecutive underpasses that contains nonzero 
many overpasses. The bridge number $b(K)$ of a knot $K$, 
is the minimum number of bridges in a knot diagram realizing $K.$ Given a Gauss diagram constructed from a 
knot diagram, the number of bridges in it is the 
number of circle arcs between two consecutive  arrow heads that contain nonzero many arrow ends. Note that 
this number is well defined even for nonrealizable Gauss diagrams, and we call it the number of bridges in a 
Gauss diagram. The {\em virtual 
bridge number\/} $\vb (K)$ of a virtual knot $K$ is the 
minimum number of bridges over all the Gauss diagrams realizing $K.$ If $\vb(K)=0,$ then $K$ is the unknot. 
Note that if $K$ is isotopic to an ordinary knot 
then $\vb(K)\leq b(K).$
It is plausible that for many $K$ we should actually have $\vb(K)=b(K),$ but we have not studied this 
question yet. 
\end{defin}

\begin{defin} An {\it isotopy\/} between two Gauss diagrams is a sequence of the Gauss diagram versions of 
Reidemeister moves. 
\end{defin}

\begin{defin}
A (generic) {\it homotopy\/} between two Gauss diagrams is a sequence of isotopies and ``flip moves'' that 
change simultaneously the direction and the sign of an arrow. If a Gauss diagram corresponds to a knot 
diagram, then the last move reflects the passage of a knot 
through  a singular
double point at a crossing identified with the arrow.
\end{defin}

Note that any two knots in $\R^3$ are homotopic as loops. However there are many different homotopy classes 
of virtual knots. In particular not every virtual knot is 
homotopic to the unknot.

\begin{defin} For an ordinary knot $K$, its {\em unknotting number\/} $u(K)$ is the minimal number of 
passages through a double point of two branches of the knot 
in a generic homotopy of $K$ to the unknot. 

For virtual knots, the unknot is the trivial Gauss diagram with no chords. If $K$ is a virtual knot homotopic 
to the unknot, then we put the {\em virtual unknotting 
number\/} $\vu(K)$ to be 
the minimal number of flip moves in a virtual homotopy that changes $K$ to the unknot. Note that if $K$ is 
nontrivial, then $\vu(K)\geq 1.$ 

Clearly if $K$ is an ordinary knot, then $\vu(K)\leq u(K).$ For many knots we have $\vu(K)=u(K).$ For 
example this is always true if $u(K)=1.$

Let $K_1, K_2$ be two different virtual knots that are homotopic to each other. We put the {\em relative 
unknotting number\/} $\rvu(K_1, K_2)$ to be the minimal 
number of flip moves 
in a generic homotopy of $K_1$ to $K_2.$ Clearly if $K_1, K_2$ are homotopic to the unknot, then we 
have $\rvu(K_1, K_2)\leq \vu(K_1)+\vu(K_2).$ Moreover if $K_2$ 
is the unknot, then $\rvu(K_1, K_2)=\vu(K_1).$ 

\end{defin}

It is well-known that there are no ordinary knots $K$ with $b(K)=1.$ However, as we show, the theory of virtual knots 
with $\vb(K)=1$ is surprisingly nontrivial.

\begin{thm}\label{main1}
For every $i\in \N$ there exists a virtual knot $K_i$ homotopic to the unknot that has  
$\vb(K_i)=1$ and $\vu(K_i)=i.$
\end{thm}

The following Theorem can be thought of a generalization of the above result to the case of virtual knots 
that are not homotopic to the unknot.

\begin{thm}\label{main2}
There are infinitely many virtual knot homotopy classes $H^{p,q}$ enumerated by pairs of 
distinct positive integers $(p,q)$ such that $H^{p,q}$ 
contains infinitely many distinct virtual knots $\{K^{p,q}_{j}\}_{j=0}^{\infty}$
with $\vb(K^{p,q}_j)=1$ and $\rvu(K^{p,q}_{j_1}, K^{p,q}_{j_2})=|j_1-j_2|$ for all $j_1, j_2.$
\end{thm}

Note that for a virtual knot homotopy class not containing the unknot, there is no obvious choice of the most 
trivial knot in the class.  
In the examples constructed in the proof of Theorem~\ref{main2} the knots $K^{p,q}_{0}$ are simpler than 
all $K^{p,q}_j, j\geq 1.$ So one can think of them as the chosen 
most trivial knots in the corresponding virtual 
homotopy classes.

\section{Proofs}

In the proofs of Theorems~\ref{main1} and~\ref{main2} we use properties of invariants studied by 
V.~Turaev~\cite{Turaev, Turaevcobordism} 
and A.~Henrich~\cite{Henrich}.

\subsection{Invariants of V.~Turaev and A.~Henrich}

\begin{defin}[Invariant $P(K)$ of A.~Henrich~\cite{Henrich}]\label{PInvariant}
Let $D$ be a Gauss diagram realizing a virtual knot $K.$ The end points of a chord $c$ in $D$ separate the 
circle in $D$ into two arcs. Choose one of the arcs and 
perform  the flip moves on all the other chords so that they point into the chosen arc. Let $i(c)$ to be the 
sum of the signs $\pm 1$ of the
chords pointing into the chosen arc after the flip moves. We put $\sign(c)=\pm 1$ to be the sign of $c.$

\def\mfoota{\footnote{Note that Goussarov, Polyak and Viro~\cite{GPV} used a different 
notion of finite order invariants of virtual knots and our invariant is not of order one in the sense of 
their work.}}

Put 
\begin{equation}
P(D)=\sum_{c \text{ such that } i(c)\neq 0} \sign(c)t^{|i(c)|} 
\end{equation}
to be the polynomial in variable $t$ 
with integer coefficients. One can show that $P(D)$ does not
depend on the choice of the diagram $D$ realizing $K$. Hence we have a polynomial  invariant $P(K)$ of a 
virtual knot $K.$ It is easy to see that $P$ is a Vassiliev-Goussarov invariant~\cite{Vassiliev},
~\cite{Goussarov1},~\cite{Goussarov2} of order one of virtual knots in the sense of 
Kauffman~\cite{Kauffman}\protect\mfoota.

This invariant, in a slightly different form, was introduced by A.~Henrich~\cite{Henrich}. It is also 
related to the cobordism invariants of knots in thickened 
surfaces studied by V.~Turaev~\cite{Turaevcobordism}.
\end{defin}

\begin{rem}[The $P$ invariant and the virtual unknotting numbers $\rvu, \vu$]\label{PRemark}
It is easy to see that a flip move on a chord $c$ with $i(c)\neq 0$ changes one of the coefficients of the 
polynomial $P$ invariant by $2$. The flip move on a chord $c$ with $i(c)=0$ does not change the $P$ invariant.
Thus, as it was observed by A.~Henrich~\cite{Henrich}, if $K_1, K_2$ are homotopic virtual knots and 
$P(K_1)-P(K_2)=\sum_{j>0 } a_j t^j$, then every 
(generic) homotopy between $K_1$ and $K_2$ involves at least $\frac{\sum_{j>0}|a_j|}{2}$ flip moves.
Thus $\rvu(K_1, K_2)\geq \frac{\sum_{j>0}|a_j|}{2}.$ 

Let $K$ be a virtual knot homotopic to the trivial knot $K_0.$ Let $P(K)=\sum_{j>0}b_jt^j$. Since $P(K_0)=0,$ 
the above relation of the $P$ invariant with the relative 
unknotting number $\rvu$ implies that the vitual unknotting number $\vu(K)\geq \frac{\sum_{j>0}|b_j|}{2}.$
\end{rem}

\begin{defin}[$u$-invariant of V.~Turaev~\cite{Turaev}] 
Let $D$ be a Gauss diagram realizing a virtual knot $K.$ Let $\overline D$ be the Gauss diagram obtained 
from $D$ by performing the flip moves that convert the signs 
of all the chords to be $+1.$ The end points of an oriented chord $c$ separate the core circle of 
$\overline D$ into two oriented arcs. The oriented arc that starts 
at the arrow head of $c$ would be called the preferred arc and denoted by $A^+_c$. The other arc is $A^-_c.$ 
Put $n_+(c)$ to be the number of oriented chords 
different from $c$ that start in $A^-_c$ and end in $A^+_c.$ Put $n_-(c)$ to be the number of oriented chords 
different from $c$ that start in $A^+_c$ and end in $A^-_c.$ Put $n(c)=n_+(c)-n_-(c)$ and put 
$\sign(n(c))=\pm 1$ to be the sign of $n(c).$

Put 
\begin{equation}
u(D)=\sum_{c\in \overline D \text{ such that } n(c)\neq 0}\sign (n(c))t^{|n(c)|}
\end{equation}
to be a polynomial in $t$ with integer coefficients. This $u(D)$ in a slightly different form was introduced 
by V.~Turaev~\cite{Turaev}. From his work one immediately gets that $u(D)$ depends 
only on the virtual homotopy class of the virtual knot $K.$ In particular, if for two virtual knots $K_1, K_2$ 
we have $u(K_1)\neq u(K_2)$, then $K_1$ and $K_2$ are not
virtually homotopic.
\end{defin}

\begin{figure}[htbp]
\begin{center}
\epsfxsize 11cm \hepsffile{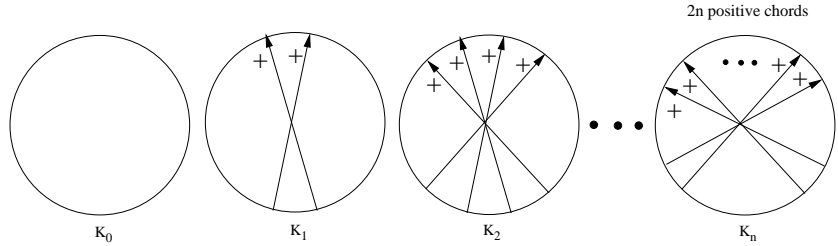}
\end{center}
\caption{Homotopic virtual bridge index one knots 
with different unknotting numbers}\label{bridgefigureone.fig}
\end{figure}

\subsection{Proof of Theorem~\ref{main1}}
Consider the virtual knots $\{K_i\}_{i=0}^{\infty}$ depicted in Figure~\ref{bridgefigureone.fig}. All these 
knots are virtually homotopic. Indeed to get from  $K_{i+1}$ to $K_i$ perform a flip move 
on the chord in $K_{i+1}$ with the rightmost arrow head. Then use the Gauss diagram version of the second 
Reidemeister move, see Figure~\ref{GaussRedemeister.fig}, 
to cancel the resulting chord and the chord with the left most arrow head in $K_{i+1}$. 
Since $K_0$ is the trivial knot, we see that $\vu(K_i)\leq i$ for all $i\in \N.$

A straightforward computation shows that $P(K_i)=2t^{2i-1}+2t^{2i-3}+\dots +2t^1.$ As it was discussed in 
Remark~\ref{PRemark} this means that $\vu(K_i)\geq \frac{2i}{2}=i.$ Hence $\vu(K_i)=i.$ 

It is easy to see that $\vb(K_i)\leq 1$ for $i\geq 1.$ Since $P(K_i)\neq 0$ for $i\geq 1$, these 
virtual knots are nontrivial and hence $\vb(K_i)=1$ for all $i\geq 1.$ 
\qed

\subsection{Proof of Theorem~\ref{main2}}
\begin{figure}[htbp]
\begin{center}
\epsfxsize 11cm \hepsffile{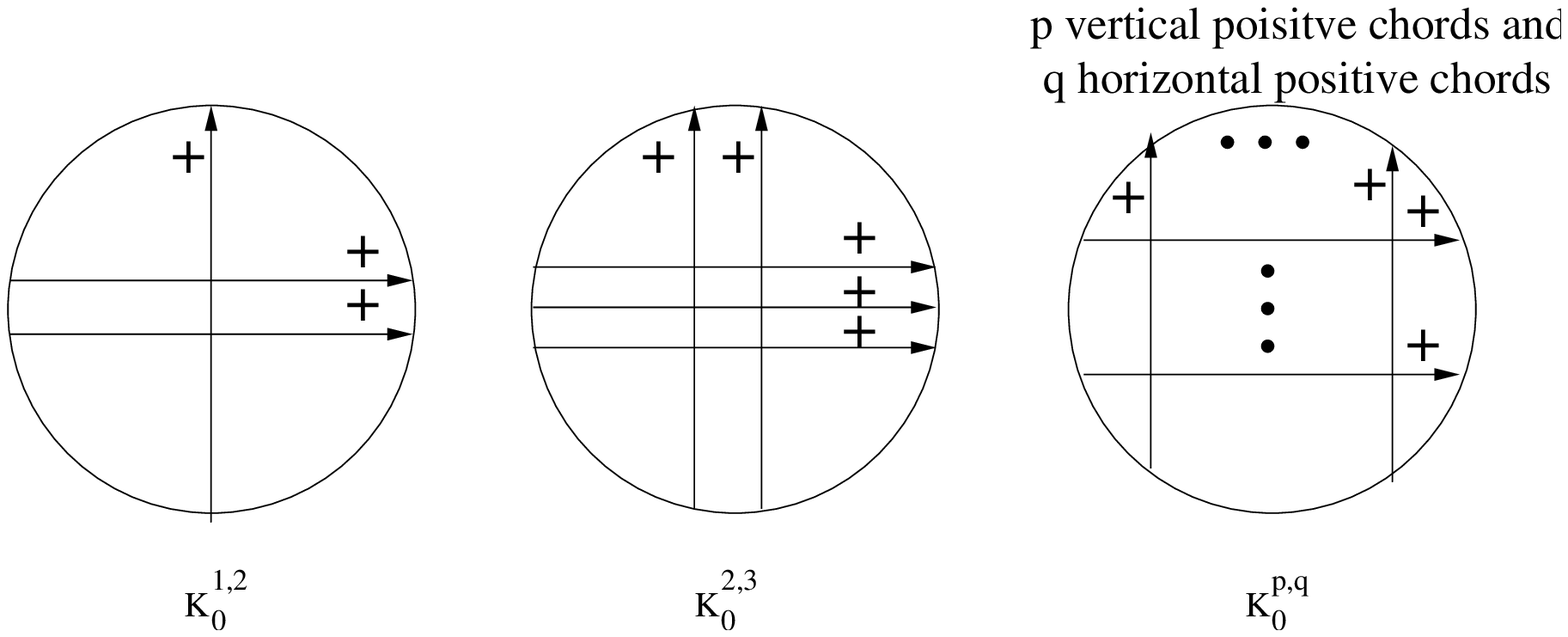}
\end{center}
\caption{nonhomotopic $K^{p,q}_0$ knots with $p\neq q$}\label{bridgefiguretwo.fig}
\end{figure}

For $p,q\in \N$ with $p\neq q$ we put $K_0^{p,q}$ to be the virtual knot shown 
in Figure~\ref{bridgefiguretwo.fig} whose Gauss diagram consists of $p$
parallel positive vertical chords oriented upwards and $q$ parallel positive horizontal chords oriented to 
the right. A straightforward computation shows that $u(K_0^{p,q})=-pt^q+qt^p.$ Since $u$ is invariant under 
virtual homotopy we see that if $p_1, q_1, p_2, q_2\in \N$ are such $p_1\neq q_1$, $p_2\neq q_2$ and 
$(p_1, q_1)\neq (p_2, q_2)$, then $K_0^{p_1, q_1}$ and $K_0^{p_2, q_2}$ are not
virtually homotopic.

\begin{figure}[htbp]
\begin{center}
\epsfxsize 11cm \hepsffile{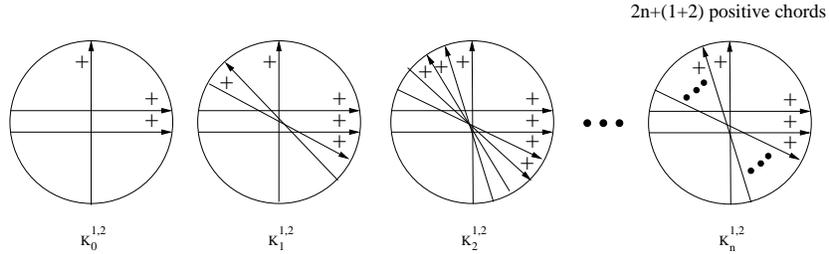}
\end{center}
\caption{Virtual bridge number one knots $K^{1,2}_n$ homotopic to $K^{1,2}_0$}\label{bridgefigurethree.fig}
\end{figure}

If we draw the horizontal and vertical arrows of $K_0^{p,q}$ to be close to the $x$- and $y$-axis, then we 
have four large arcs of 
the circle in $K^{p,q}_0$ that are free of chord ends. We will refer to these arcs as the arcs located in 
the corresponding four quadrants. For $n\in \N$ put $K_n^{p,q}$ to be the virtual knot whose Gauss diagram 
is obtained from the diagram of $K_0^{p,q}$ by adding $2n$ positive chords that pass through the circle 
center. The first $n$ of the added chords should start in the fourth quadrant arc, end in the second quadrant
arc and have almost vertical slope.  The other $n$ chords should start in the second quadrant arc, end in 
the fourth quadrant arc and have almost horizontal slope, see Figure~\ref{bridgefigurethree.fig}. 

\begin{figure}[htbp]
\begin{center}
\epsfxsize 12cm \hepsffile{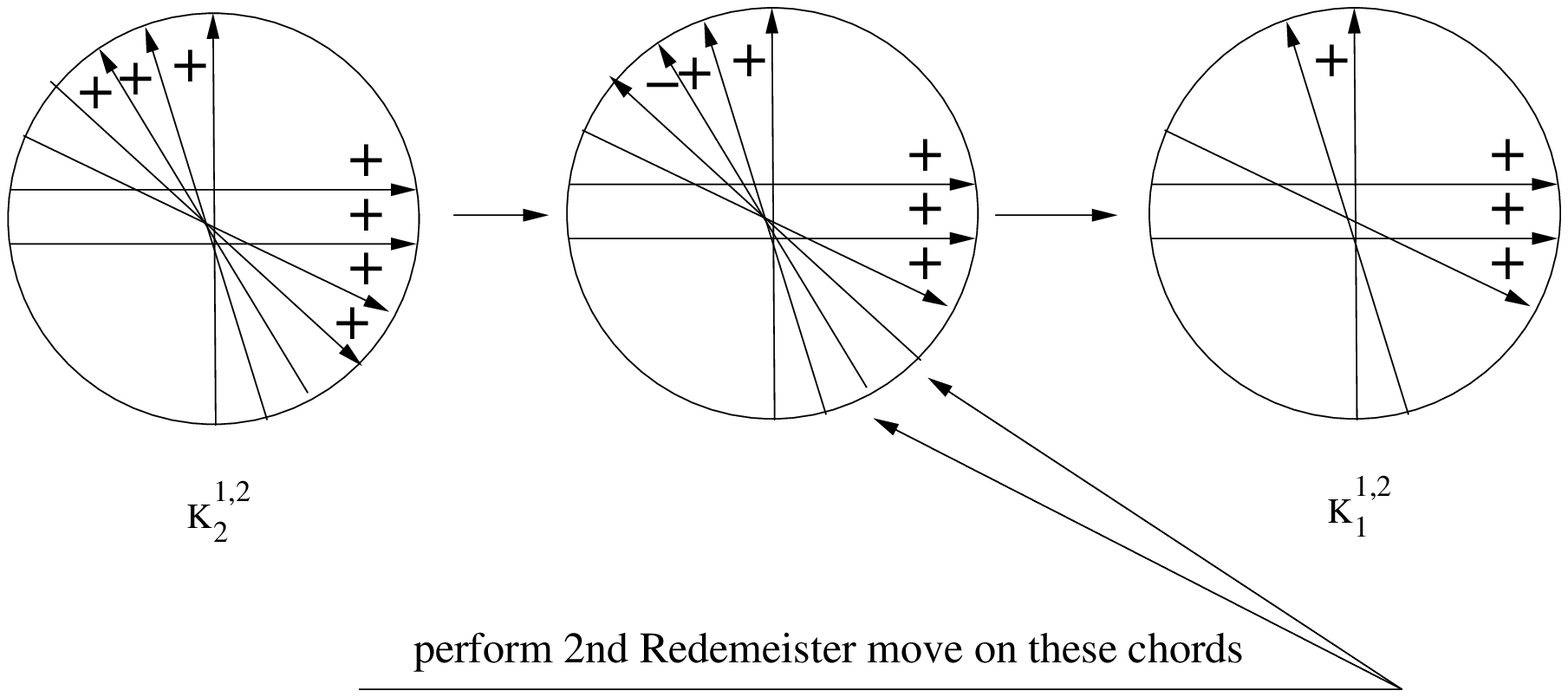}
\end{center}
\caption{Homotopy between $K^{1,2}_2$ and $K^{1,2}_1$ knots}\label{bridgefigurefour.fig}
\end{figure}

We claim that for fixed $p,q$ all the knots $K_n^{p,q}$ are virtually homotopic. For $n\geq 1$ the homotopy
between $K_{n}^{p,q}$ and $K_{n-1}^{p,q}$ is constructed as follows. Take the chord 
from the group of $n$ chords going from the second to the fourth quadrant arc whose slope is the farthest 
from 
the horizontal one and perform the flip move on it. Use the Gauss diagram version of the second 
Reidemiester move to cancel the resulting chord and the chord in the group of $n$ chords going from the 
fourth into the second quadrant arc whose slop is the farthest from the vertical, 
see Figure~\ref{bridgefigurefour.fig}. Since this virtual homotopy between $K_{n}^{p,q}$ and $K_{n-1}^{p,q}$ 
involves only one flip move,  we have 
$\rvu(K_{n_1}^{p,q}, K_{n_2}^{p,q})\leq |n_1-n_2|,$ for all $n_1, n_2$.

A straightforward computation shows that 
\begin{equation}
\begin{split}
P(K_{n}^{p,q})=\\
t^{2n+p+q-1}+\dots+t^{2n+p+q-(2n-1)}+
qt^p+pt^q+t^{p+q+1}+\dots t^{p+q+(2n-1)}=\\
qt^p+pt^q+\sum_{i=1}^n2t^{p+q+(2i-1)}.
\end{split}
\end{equation}

Thus 
\begin{equation}
P(K_{n_1}^{p,q})-P(K_{n_2}^{p,q})=\sum_{i=\min(n_1, n_2)+1}^{i=\max(n_1, n_2)}2t^{p+q+2i-1}.
\end{equation}
Hence by Remark~\ref{PRemark}, $\rvu(K_{n_1}^{p,q}, K_{n_2}^{p,q})\geq \frac{2|n_2-n_1|}{2}=|n_2-n_1|.$
Thus we have $\rvu(K_{n_1}^{p,q}, K_{n_2}^{p,q})=|n_2-n_1|.$ 

It is easy to see that $\vb(K^{p,q}_n)\leq 1.$ Since $P(K^{p,q}_n)\neq 0$, these knots are nontrivial 
and we have 
$\vb(K^{p,q}_n)=1.$
\qed.

{\bf Acknowledgments.} 
The second author is very thankful to Alina Vdovina for the very useful mathematical discussions 
during which the notion of the virtual bridge number was introduced.

\end{document}